\documentclass[12pt]{amsart}
\usepackage{a4wide,amsmath,amssymb,graphicx}
\usepackage{color}

\let\pa\partial
\let\na\nabla
\let\eps\varepsilon
\newcommand{\R}{\mathbb{R}}

\newcommand{\diver}{\textnormal{div}}

\newtheorem{theorem}{Theorem}

\theoremstyle{remark}

\begin{document}
\title[A simplified energy-transport model]{Existence analysis for a 
simplified transient energy-transport model for semiconductors}

\author{Ansgar J\"ungel}
\address{Institute for Analysis and Scientific Computing, Vienna University of  
	Technology, Wiedner Hauptstra\ss e 8--10, 1040 Wien, Austria}
\email{juengel@tuwien.ac.at} 

\author{Ren\'e Pinnau}
\address{Fachbereich Mathematik, Technische Universit\"at Kaiserslautern,
Erwin-Schr\"odinger-Stra\ss e, 67663 Kaiserslautern, Germany}
\email{pinnau@mathematik.uni-kl.de}

\author{Elisa R\"ohrig}
\address{ITWM, Fraunhofer-Zentrum, Fraunhofer-Platz 1, 67663 Kaiserslautern, Germany}
\email{elisa.roehrig@itwm.fraunhofer.de}

\date{\today}

\thanks{The first author acknowledges partial support from   
the Austrian Science Fund (FWF), grants P20214, P22108, and I395, and    
the Austrian-French Project of the Austrian Exchange Service (\"OAD).
The second author acknowledges support from the German Science Foundation (DFG), 
grants PI 408/5 and PI 408/7 in the context of the SPP 1253.} 

\begin{abstract}
A simplified transient energy-transport system for semiconductors subject to mixed
Dirichlet-Neumann boundary conditions is analyzed. 
The model is formally derived from the non-isothermal hydrodynamic equations in a
particular vanishing momentum relaxation limit. It consists of a drift-diffusion-type
equation for the electron density, involving temperature gradients,
a nonlinear heat equation for the electron temperature, and the Poisson equation
for the electric potential. The global-in-time existence of bounded weak solutions
is proved. The proof is based on the Stampacchia
truncation method and a careful use of the temperature equation.
Under some regularity assumptions on the gradients of the variables,
the uniqueness of solutions is shown. Finally, numerical simulations for a ballistic
diode in one space dimension illustrate the behavior of the solutions.
\end{abstract}

 

\maketitle


\section{Introduction}

The basic model for the charge transport in semiconductor devices
are the drift-diffusion equations for the electron density and the electric potential.
This model gives fast and satisfactory simulation results for devices
on the micrometer scale, but it is not able to cope with so-called hot-electron
effects in nanoscale devices. A possible solution is to incorporate the mean
energy in the model equations, which leads to energy-transport equations,
first presented by Stratton \cite{Str62} and later derived from the semiconductor
Boltzmann equation by Ben Abdallah and Degond \cite{BDG96}. The analysis of the
energy-transport model 
is very involved due to the strong coupling and temperature gradients.
Therefore, we consider in this paper a simplified energy-transport model
which still includes temperature gradients but the coupling to the energy
equation is weaker than in the full model. An important feature of our model
is that it is derived formally from the hydrodynamic semiconductor equations
in a zero relaxation time limit, which provides a physical modeling basis without
heuristics (see Section \ref{sec.deriv}). Our goal is to prove the
existence and uniqueness of solutions to this model and to provide some
numerical illustrations of the solutions.

The model consists of a drift-diffusion-type equation
for the electron density $n(x,t)$, a nonlinear heat
equation for the electron temperature $\theta(x,t)$, and the Poisson equation
for the electric potential $V(x,t)$:
\begin{align}
  \pa_t n -\diver(\na(n\theta)+n\na V) &= 0, \label{1.eq.n} \\
  \diver(\kappa(n)\na\theta) &= \frac{n}{\tau}(\theta-\theta_L(x)), 
  \label{1.eq.t}\\
  -\lambda^2\Delta V &= n-C(x) \quad\mbox{in }\Omega,\ t>0. \label{1.eq.v}
\end{align}
Here, $\kappa(n)$ is the thermal conductivity, $\theta_L(x)$ the lattice
temperature, and $C(x)$ the doping profile characterizing the device under
consideration. The scaled physical parameters are the energy relaxation time
$\tau>0$ and the Debye length $\lambda>0$. Equations \eqref{1.eq.n}-\eqref{1.eq.v}
hold in the bounded domain $\Omega\subset\R^d$ ($d\ge 1$) with the initial condition
\begin{equation}\label{1.ic}
  n(0) = n_I \quad\mbox{in }\Omega.
\end{equation}
We suppose that the boundary $\pa\Omega\in C^{0,1}$ consists of two
parts $\Gamma_D$ and $\Gamma_N$ satisfying $\pa\Omega=\Gamma_D\cup\Gamma_N$,
$\Gamma_D\cap\Gamma_N=\emptyset$, $\Gamma_N$ is closed, and the $(d-1)$-dimensional
Lebesgue measure of $\Gamma_D$ is positive.
The electron density, temperature, and potential are assumed to be known on the
Dirichlet boundary, which models the contacts, 
whereas the Neumann boundary models insulated boundary parts:
\begin{equation}\label{1.bc}
  \begin{aligned}
  n=n_D, \quad \theta=\theta_D, \quad V=V_D \quad&\mbox{on }\Gamma_D, \\
  \na n\cdot\nu = \na\theta\cdot\nu = \na V\cdot\nu = 0 \quad&\mbox{on }\Gamma_N, 
  \end{aligned}
\end{equation}
where $\nu$ denotes the exterior unit normal vector on $\pa\Omega$. 

Before we detail our analytical results, we review related models in the literature.
First, temperature effects have been included in the drift-diffusion equations 
by allowing for temperature-dependent diffusivities \cite{SeTr85} or
temperature-dependent mobilities \cite{FaWu01,GuWu07,WuXu06,Yin95}
coupled to a heat equation. Typically, 
the so-called non-isothermal drift-diffusion equations are of the form
\begin{align*}
  & \pa_t n - \diver\, J_n = 0, \quad J_n = D\na n + \mu n\na V, \\
  & \pa_t\theta - \diver(\kappa(\theta)\na\theta) = F, \quad
  F = J_n\cdot\na V + W,
\end{align*}
where $J_n$ is the particle current density, $D$ and $\mu$ are the diffusivity and
mobility, respectively, and $W=-n(\theta-\theta_L)/\tau$ is the relaxation term.
The difficulty in these models is that the Joule heating term $J_n\cdot\na V$ involves
quadratic gradients of the potential, which resembles the thermistor problem;
see, e.g., \cite{Xu93}. However, temperature gradients in $J_n$, which need to
be taken into account, have been ignored. 

In \cite{Xu09}, Xu allowed for temperature gradients in $J_n$ but he
truncated, as in \cite{Yin95}, the Joule heating term by setting 
$F=\max\{0,J_n\cdot\na V+W\}$ in order to allow for a maximum principle.
A different approach was adopted in \cite{Bar98}, where 
a kind of quasi-Fermi potential via $\phi = n\exp(-V/\theta)$ was introduced.
Ths model of \cite{Bar98} 
includes temperature gradients, but the coefficient contains the
electric potential which is not the case in the energy-transport models
derived in \cite{BDG96}.
We also mention non-isothermal systems with simplified thermodynamic forces
which were studied in \cite{AlXi94}. 

Compared to our model \eqref{1.eq.n}-\eqref{1.eq.v},
the energy-transport equations contain cross-diffusion terms also in the 
energy equation \cite{Jue10}. A typical form of these models reads as
\begin{align*}
  & \pa_t n - \diver\, J_n = 0, \quad 
  J_n = \na(n\theta^\alpha) + n\theta^{\alpha-1}\na V,\\
  & \frac32\pa_t(n\theta) - \diver J_e = J_n\cdot\na V + W, 
  \quad J_e = \na(n\theta^{\alpha+1}) + n\theta^\alpha\na V,
\end{align*}
where the parameter $\alpha>0$ is related to the elastic scattering rate
in the collision operator (see Example 6.8 in \cite{Jue09}). In our model
\eqref{1.eq.n}-\eqref{1.eq.t}, $\alpha=1$, and the diffusion scaling implies 
that the variation of the energy density, $\frac32\pa_t(n\theta)$, and the Joule
heating term are negligible (see Section \ref{sec.deriv}). 
The main difficulty of the above model is that the corresponding diffusion matrix
is neither diagonal nor tridiagonal and that it degenerates for $n=0$ or $\theta=0$.
Existence results were achieved for stationary equations near thermal equilibrium
\cite{FaIt01,Gri99} and for the transient model \cite{ChHs03,CHL04,Yon09} 
if the initial data are close to the stationary drift-diffusion solutions.
General existence results, both for the stationary and time-dependent model,
were proved in \cite{DGJ97,DGJ98} but the diffusion matrix was assumed
to be uniformly positive definite, thus avoiding the degeneracy.
All these results give only partial answers to the well-posedness of the 
problem, and a complete global existence theory for the energy-transport equations
for any data and with physical transport coefficients is still missing.

In this paper, we wish to bring forward the existence theory for
energy-transport-type models by analyzing the
system \eqref{1.eq.n}-\eqref{1.eq.v}, whose complexity is in between the 
well-understood drift-diffusion model and the energy-transport equations. 
In fact, in our model, the energy equation simplifies such that the application of the
maximum principle for $\theta$ becomes possible. The remaining difficulties
are due to the drift term $n\na\theta$ in \eqref{1.eq.n} and the
quasilinearity $\kappa(n)$ in \eqref{1.eq.t}. Note that, in view of
the mixed boundary conditions, we cannot expect the regularity 
$\na\theta\in L^\infty$ which would simplify the existence proofs significantly.

Our main idea is a careful use of the temperature equation in order to
deal with the drift term $n\na\theta$. More precisely, 
we replace this term in \eqref{1.eq.n} formally by
$$
  \diver(n\na\theta) 
  = \diver\left(\frac{n}{\kappa}\kappa\na\theta\right)
  = \frac{n}{\kappa}\diver(\kappa\na\theta)
  + \na n\cdot\na\theta - \frac{n}{\kappa}\na\theta\cdot
  \left(\frac{\pa\kappa}{\pa n}\na n + \frac{\pa\kappa}{\pa\theta}\na\theta\right),
$$
and using \eqref{1.eq.t}, we find that \eqref{1.eq.n} equals
\begin{equation}\label{1.n2}
  \pa_t n - \diver(\theta\na n) 
  = \diver(n\na V) + \frac{n^2}{\kappa}(\theta-\theta_L) 
  + \left(1 - \frac{n}{\kappa}\frac{\pa\kappa}{\pa n}\right)\na n\cdot\na\theta
  - \frac{n}{\kappa}\frac{\pa\kappa}{\pa\theta}|\na\theta|^2.
\end{equation}
The computations will be made rigorous on a weak formulation level
in Section \ref{sec.ex}. From the above formulation we see that the last term
on the right-hand side models a sink if $\pa\kappa/\pa\theta\ge 0$. 
This condition is satisfied, for instance, in the case of the Wiedemann-Franz model
$\kappa(n,\theta)=n\theta$. By the maximum principle, we expect to obtain
an upper bound for $n$.

However, we need the stronger condition $\pa\kappa/\pa\theta=0$.
The reason is that the lack of time regularity for $\theta$ makes it difficult
to deal with nonlinear terms, such as $\theta\na n$, 
to prove the continuity of the fixed-point operator.
Although in physical models, it is often assumed that the thermal
conductivity depends on the temperature $\theta$, a dependency on $n$ only
also occurs in the physical literature. 
For instance, the choice $\kappa(n)=n$ was suggested in 
\cite[Formula (2.16)]{JeSh95} to study
spurious velocity overshoots in hydrodynamic semiconductor models.

From the physical application, we expect that the electron density $n$ stays
positive if it is positive initially and on the Dirichlet boundary parts.
Even if $\kappa$ depends on $n$ only, the proof of a positive lower bound for $n$
is not obvious, since it is not clear how to deal with the term 
$\na n\cdot\na\theta$ in \eqref{1.n2} which is in $L^1$ only. 
We suppose that either $\kappa(n)$ is strictly positive or $\kappa(n)=n$. 
In the former case, we avoid any degeneracy; in the latter case, 
$(n/\kappa)(\pa\kappa/\pa n)=1$, 
and the term involving $\na n\cdot\na\theta$ in \eqref{1.n2} vanishes.

Motivated by the above considerations, we impose the following conditions 
on the thermal conductivity: Let $\kappa\in C^1([0,\infty))$ such that
there exist $\kappa_0$, $\kappa_1$, $n_*$, $n^*>0$ with
\begin{equation}\label{1.kappa}
  \begin{aligned}
  & \mbox{(i) }\kappa(z)>0\mbox{ for all }z>0, \\
  & \mbox{(ii) either }\kappa(z)\ge\kappa_0>0\mbox{ for all }z\ge 0, \mbox{ or }
  \kappa(z)=z\mbox{ for all }0\le z\le n_*; \\
  & \mbox{(iii) }\kappa(z)\ge \kappa_1 z\mbox{ for all }z\ge n^*.
  \end{aligned}
\end{equation}
Condition (i) allows for the degenerate case $\kappa(0)=0$.
Condition (ii) ensures the uniform ellipticity of equation
\eqref{1.eq.t}. Indeed, if $\kappa(n)=n$ for $n\le n_*$, we are able to prove that the
solution $n$ is strictly positive and then, $\kappa(n)$ is strictly positive, too.
The last condition is needed to prove an upper bound for the particle density. 

The boundary data are assumed to satisfy
\begin{equation}\label{1.bound}
  \begin{aligned}
  & n_D,\ V_D\in L^2(0,T;H^1(\Omega)), \quad \theta_D\in L^q(0,T;W^{1,q}(\Omega)), \\
  & n_D,\ \theta_D\in L^\infty(0,T;L^\infty(\Omega)), \quad
  \inf_{\Omega_T}n_D > 0,\ \inf_{\Omega_T}\theta_D>0,
  \end{aligned}
\end{equation}
where $q>2$ and $\Omega_T=\Omega\times(0,T)$. 
The initial data and the given functions fulfill the conditions
\begin{equation}\label{1.init}
  n_I,\,\theta_L,\,C\in L^\infty(\Omega), \quad \inf_\Omega n_I>0, \quad
  \inf_\Omega \theta_L>0, \quad \inf_\Omega C(x)\ge 0.
\end{equation}
In order to deal with the mixed Dirichlet-Neumann conditions,
we introduce the space
$$
  H_0^1(\Omega\cup\Gamma_N) = \{u\in H^1(\Omega): u = 0\mbox{ on }\Gamma_D\}.
$$
For properties of this space, we refer to \cite[Chapter 1.7.2]{Tro87}.
Furthermore, we set $H^{-1}(\Omega\cup\Gamma_N) = (H_0^1(\Omega\cup\Gamma_N))'$.

\begin{theorem}[Existence of solutions]\label{thm.ex1}
Let $\Omega\subset\R^d$ ($d\ge 1$) be a bounded domain with $\pa\Omega\in C^{0,1}$,
$T$, $\tau$, $\lambda>0$, and let
$\kappa\in C^1([0,\infty))$ satisfy \eqref{1.kappa}.
Furthermore, assume that \eqref{1.bound} and \eqref{1.init} hold.
Then there exists a weak solution $(n,\theta,V)\in L^2(0,T;H^1(\Omega))^3$
to \eqref{1.eq.n}-\eqref{1.bc} satisfying
$\pa_t n\in L^2(0,T;H^{-1}(\Omega\cup\Gamma_N))$ and
$$
  0 \le n(t)\le K_0 e^{\beta t}, \quad 
  0 < m \le \theta(t) \le M \quad\mbox{in }\Omega,\ t\in(0,T).
$$
Furthermore, if $\kappa(z)=z$ for $0\le z\le n_*$,
$$
  n(t) \ge k_0 e^{-\alpha t} > 0 \quad\mbox{in }\Omega,\ t\in(0,T).
$$
\end{theorem}

In the above theorem, the constants are defined by
\begin{equation}\label{1.const}
\begin{aligned}
  K_0 &= \max\left\{n^*,\ \sup_\Omega n_I, 
  \sup_{\Gamma_D\times(0,T)}n_D,\ \sup_\Omega C(x)\right\}, \\
  k_0 &= \min\left\{n_*,\ \inf_\Omega n_I,\ \inf_{\Gamma_D\times(0,T)}n_D\right\}, \\
  M &= \max\left\{\sup_\Omega \theta_L, \sup_{\Gamma_D\times(0,T)}\theta_D\right\},
  \quad m = \min\left\{\inf_\Omega \theta_L,\ 
  \inf_{\Gamma_D\times(0,T)}\theta_D\right\}, \\
  \alpha &= \frac{1}{\tau}\sup_\Omega\theta_L + \frac{1}{\lambda^2}, 
  \quad \beta = \frac{M}{\tau \kappa_1}.
  \end{aligned}
\end{equation}
The proof of the theorem is based on the Leray-Schauder fixed-point theorem
and the Stampacchia truncation method. In particular, the truncation is needed in the
diffusion coefficients of $\diver(\theta\na n)$ and $\diver(\kappa(n)\na\theta)$
to make these expressions uniformly elliptic.

Due to the quasilinearity of the temperature equation, we are able to show 
the uniqueness of solutions only in a function space which includes
bounded gradients.

\begin{theorem}[Uniqueness of solutions]\label{thm.unique}
Let the assumptions of Theorem \ref{thm.ex1} hold and let 
$\kappa$ be locally Lipschitz continuous on $[0,\infty)$.
Then there exists a unique solution $(n,\theta,V)$
to \eqref{1.eq.n}-\eqref{1.bc} in the class of bounded weak solutions
satisfying $n\in L^\infty(0,T;W^{1,p}(\Omega))$, 
$\theta\in L^\infty(0,T;W^{1,\infty}(\Omega))$, where $p>2$ if $d=2$ and
$p\ge d$ if $d\ge 3$.
\end{theorem}

The paper is organized as follows. Equations \eqref{1.eq.n}-\eqref{1.eq.v}
are formally derived from the hydrodynamic model in Section \ref{sec.deriv}.
The existence theorem is proved in Section \ref{sec.ex},
and Section \ref{sec.unique} is devoted to the proof of the uniqueness theorem.
In Section \ref{sec.num}, we present numerical results for a simple
one-dimensional ballistic diode illustrating the behavior of the electron temperature
in the presence of a cooling and heating lattice temperature.


\section{Derivation of the model equations}\label{sec.deriv}

Equations \eqref{1.eq.n}-\eqref{1.eq.v} are formally derived from
the (scaled) hydrodynamic model (see, e.g., \cite[Chapter 9]{Jue09}):
\begin{align*}
  \pa_t n - \diver\, J_n &= 0, \\
  \pa_t J_n - \diver\left(\frac{J_n\otimes J_n}{n}\right) - \na(n\theta) - n\na V
  &= -\frac{J_n}{\tau_p}, \\
  \pa_t(ne) - \diver(J_n(e+\theta)) - J_n\cdot\na V - \diver(\kappa(n,\theta)\na\theta)
  &= -\frac{n}{\tau_e}\left(e-\frac32 \theta_L\right), 
\end{align*}
and $V$ is given by the Poisson equation \eqref{1.eq.v}.
Here, $J_n$ denotes the particle current density, $J_n\otimes J_n$ is a tensor product, 
$\tau_p$ is the momentum relaxation time, and
$\tau_e$ the energy relaxation time.
The energy density is the sum of the thermal and kinetic energies:
$$
  ne = \frac32 n\theta + \frac{|J_n|^2}{2n}.
$$

Energy-transport equations can be derived from the vanishing momentum
relaxation limit. To this end, we set $\eps=\tau_p$ and rescale the equations by 
$t\to t/\eps$ and $J\to \eps J$. This corresponds to the physical situation
of a long time scale and small current densities. 
The rescaled equations become:
\begin{align}
  & \pa_t n - \diver\, J_n = 0, \quad 
  ne = \frac32 n\theta + \frac {\eps^2}{2}\,\frac{|J_n|^2}{n}, \label{2.n} \\
  & \eps^2\pa_t J_n - \eps^2\diver\left(\frac{J_n\otimes J_n}{n}\right) 
  - \na(n\theta) - n\na V = -J_n, \label{2.j} \\
  & \eps\pa_t(ne) - \eps\diver(J_n(e+\theta)) - \eps J_n\cdot\na V 
  - \diver(\kappa(n,\theta)\na\theta)
  = -\frac{n}{\tau_e}\left(e-\frac32 \theta_L\right). \nonumber
\end{align}
In the formal limit $\eps\to 0$, we obtain the limiting model 
$$
  \pa_t n - \diver(\na(n\theta)+n\na V), \quad
  \diver(\kappa(n,\theta)\na\theta)
  = \frac{n}{\tau}\left(e-\frac32\theta_L\right), \quad e = \frac32 \theta,
$$
which corresponds to \eqref{1.eq.n}-\eqref{1.eq.t} with $\tau=2\tau_e/3$.

In the literature, usually a different limit is performed in order to
derive energy-transport equations. Indeed,
if we rescale additionally $\kappa\to \eps\kappa$ (small thermal conductivity),
and assume that the energy relaxation time is of the same order as the
momentum relaxation time, $\tau=\tau_0=\eps$, the rescaled energy equation
reads as
\begin{equation}\label{2.ne.limit}
  \eps\pa_t(ne) - \eps\diver(J_n(e+\theta)) - \eps J_n\cdot\na V 
  - \eps\diver(\kappa(n,\theta)\na\theta)
  = -\eps n\left(e-\frac32 \theta\right).
\end{equation}
Then, dividing this equation by $\eps$ and 
performing the formal limit $\eps\to 0$ in \eqref{2.n} and \eqref{2.j}, we find 
the usual energy-transport model 
with particular diffusion coefficients (see \cite[Chapter 6.4]{Jue09}).

Our simplified model is valid in diffusive situations in which the thermal
conductivity is strong and the energy relaxation time is much larger than
the momentum relaxation time. 


\section{Proof of Theorem \ref{thm.ex1}}\label{sec.ex}

The existence proof is based on the Leray-Schauder fixed-point
theorem and a truncation method. For this, we consider the truncated problem
\begin{align}
  \pa_t n -\diver(\theta_{m,M}\na n) &= \diver(n_K\na(\theta+V)), \label{3.eq.n} \\
  \diver(\kappa(n_{k,K})\na\theta) 
  &= \frac{n_K}{\tau}(\theta-\theta_L), \label{3.eq.t}\\
  -\lambda^2\Delta V &= n_K-C(x) \quad\mbox{in }\Omega,\ t>0, \label{3.eq.v}
\end{align}
with the initial and boundary conditions \eqref{1.ic}-\eqref{1.bc},
where 
\begin{align*}
  & n_K = \max\big\{0,\min\{K,n\}\big\}, \\
  & n_{k,K} = \max\big\{k,\min\{K,n\}\big\}, \\
  & \theta_{m,M} = \max\big\{m,\min\{M,\theta\}\big\},
\end{align*}
and $k=k(t)=k_0 e^{-\alpha t}$, $K=K(t)=K_0 e^{\beta t}$.
We recall that the constants $m$, $M$, $k_0$, $K_0$, $\alpha$, and $\beta$ 
are defined in \eqref{1.const}. Observe that the lower truncation of $n$ in
\eqref{3.eq.t} is not necessary if $\kappa(n)\ge\kappa_0>0$ for all $n\ge 0$.
In this case, we replace $\kappa(n_{k,K})$ by $\kappa(n_K)$.

We divide the proof in several steps.

{\em Step 1: Definition of the fixed-point operator.}
Let $w\in L^2(0,T;L^2(\Omega))$ and $\sigma\in[0,1]$. For given $t\in(0,T)$,
let $V(t)\in H^1(\Omega)$ be the unique solution to the linear problem
$$
  -\lambda^2\Delta V(t) = w(t)_K-C(x)\mbox{ in }\Omega, \quad
  V(t)=V_D(t)\mbox{ on }\Gamma_D, \quad 
  \na V(t)\cdot\nu=0\mbox{ on }\Gamma_N.
$$
Since $w\in L^2(0,T;L^2(\Omega))$, we find that $V:(0,T)\to H^1(\Omega)$
is Bochner-measurable and $V\in L^2(0,T;H^1(\Omega))$ (see, e.g., 
\cite[pp.~1133f.]{AnCh94}).

Next, let $\theta(t)\in H^1(\Omega)$ be the unique solution to the
linear uniformly elliptic problem
$$
  \diver\big(\kappa(w(t)_{k,K})\na\theta\big) 
  = \frac{w(t)_K}{\tau}(\theta-\theta_L)\mbox{ in }\Omega, \quad
  \theta=\theta_D(t)\mbox{ on }\Gamma_D, \quad 
  \na \theta\cdot\nu=0\mbox{ on }\Gamma_N.
$$
Again, the integrability of $w$ allows us to conclude that 
$\theta\in L^2(0,T;H^1(\Omega))$. 

Finally, consider the linear parabolic problem
\begin{align*}
  & \pa_t n - \diver(\theta_{m,M}\na n) = \sigma\diver\big(w_K\na(\theta+V)\big)
  \quad\mbox{in }\Omega,\ t>0, \\
  & n=\sigma n_D\mbox{ on }\Gamma_D, \quad \na n\cdot\nu = 0\mbox{ on }\Gamma_N,
  \quad n(0)=\sigma n_I\mbox{ in }\Omega.
\end{align*}
Since the right-hand side of the parabolic equation is an element of
$L^2(0,T;H^{-1}(\Omega\cup\Gamma_N))$, there exists a unique solution
$n\in L^2(0,T;H^1(\Omega))\cap H^1(0,T;H^{-1}(\Omega\cup\Gamma_N))$.
This shows that the operator $S:L^2(0,T;L^2(\Omega))\times[0,1]\to
L^2(0,T;L^2(\Omega))$, $(w,\sigma)\mapsto n$, is well defined.
It holds that $S(w,0)=0$ for all $w\in L^2(0,T;L^2(\Omega))$.

By using $\theta-\theta_D$ as a test function in \eqref{3.eq.t}, standard
estimates and the lower bound of $\kappa$ show that
$$
  \|\theta\|_{L^2(0,T;H^1(\Omega))} \le c_1,
$$
where $c_1>0$ depends on $\kappa_0$, $m$, $M$, $K$, $\theta_L$, and $\theta_D$.
Similarly, 
$$
  \|V\|_{L^2(0,T;H^1(\Omega))} \le c_2,
$$
where $c_2>0$ depends on $K$, $\lambda$, $C(x)$, and $V_D$. Therefore,
employing $n-\sigma n_D$ as a test function in \eqref{1.eq.n}, a Gronwall
estimate implies that
$$
  \|n\|_{L^2(0,T;H^1(\Omega))} 
  + \|\pa_t n\|_{L^2(0,T;H^{-1}(\Omega\cap\Gamma_N))} \le c_3,
$$
where $c_3>0$ depends on $m$, $K$, $n_D$, $n_I$, $c_1$, and $c_2$.

We claim that $\theta$ is slightly more regular. Indeed,
using the (admissible) test function $(\theta-M)^+=\max\{M,\theta\}$ in \eqref{3.eq.t},
we obtain
\begin{align*}
  \kappa_*\int_\Omega|\na(\theta-M)^+|^2 dx
  &\le \int_\Omega\kappa(n_{k,K})|\na(\theta-M)^+|^2 dx \\
  &= -\frac{1}{\tau}\int_\Omega n_K(\theta-\theta_L)(\theta-M)^+ \le 0,
\end{align*}
since $\theta-\theta_L\ge 0$ on $\{\theta>M\}$,
where $\kappa_*=\kappa_0>0$ or $\kappa_*=\min_{z\in[k,K]}\kappa(z)>0$ 
(see \eqref{1.kappa}). We infer that $\theta\le M$
on $\Omega$, $t>0$. In a similar way, the test function 
$(\theta-m)^-=\min\{m,\theta\}$ yields $\theta\ge m$. In particular,
we have $\theta_{m,M}=\theta$. Thus, the right-hand side of the heat equation
is an element of $L^\infty(0,T;L^\infty(\Omega))$. By elliptic regularity, 
we have \cite[Theorem 1]{Gro89} $\theta(t)\in W^{1,p}(\Omega)$ for some
$2<p\le q$, and hence, $\na\theta\in L^p(0,T;L^p(\Omega))$.

{\em Step 2: Continuity of the fixed-point operator.}
Let $w_j\to w$ strongly in $L^2(0,T;L^2(\Omega))$ and $\sigma_j\to \sigma$
as $j\to\infty$. Let $\theta_j$ and $V_j$ be the solutions to 
\begin{equation}\label{3.eq.j}
  \diver(\kappa((w_j)_{k,K})\na\theta_j) = \frac{(w_j)_K}{\tau}(\theta_j-\theta_L),
  \quad -\lambda^2\Delta V_j = (w_j)_K-C(x) \quad\mbox{in }\Omega,
\end{equation}
with the corresponding boundary conditions.
Then, by the above elliptic estimates, up to a subsequence,
$$
  \theta_j\rightharpoonup \theta, \quad V_j\rightharpoonup V
  \quad\mbox{weakly in }L^2(0,T;H^1(\Omega)).
$$
Since $\kappa((w_j)_{k,K})\to\kappa(w_{k,K})$ strongly in $L^r(0,T;L^r(\Omega))$ for
any $r<\infty$, we can pass to the limit in \eqref{3.eq.j} to obtain
$$
  \diver(\kappa(w_{k,K})\na\theta) = \frac{w_K}{\tau}(\theta-\theta_L),
  \quad -\lambda^2\Delta V = w_K-C(x) \quad\mbox{in }\Omega.
$$

In view of the compact embedding $H^1(\Omega)\hookrightarrow L^2(\Omega)$,
Aubin's lemma shows that 
$L^2(0,T;$ $H^1(\Omega))\cap H^1(0,T;H^{-1}(\Omega\cup\Gamma_N))$ 
is compactly embedded into $L^2(0,T;L^2(\Omega))$. Thus, the above estimate for $n_j$
proves that, again up to a subsequence,
$$
  n_j\to n\quad\mbox{strongly in }L^2(0,T;L^2(\Omega)).
$$
We have to show that $n=S(w,\sigma)$. This is proved by passing to the
limit in the parabolic equation satisfied by $n_j$. The problem is the limit
of $(\theta_j\na n_j)$ since $(\theta_j)$ and $(\na n_j)$ both
converge only weakly. We claim that $\theta_j\to\theta$ strongly
in $L^2(0,T;H^1(\Omega))$. Taking the difference of the equations satisfied
by $\theta_j$ and $\theta$, respectively, and using $\theta_j-\theta$ as
a test function, we find that 
\begin{align*}
  \int_0^T\int_\Omega & \kappa((w_j)_{k,K})|\na(\theta_j-\theta)|^2 dxdt \\
  &= -\int_0^T\int_\Omega\big(\kappa((w_j)_{k,K})-\kappa(w_{k,K})\big)\na\theta
  \cdot\na(\theta_j-\theta)dxdt \\
  &\phantom{xx}{}- \frac{1}{\tau}\int_0^T\int_\Omega\big(((w_j)_K-w_K)\theta 
  + (w_j)_K(\theta_j-\theta) - ((w_j)_K-w_K)\theta_L\big)(\theta_j-\theta)dxdt
\end{align*}

\begin{align*}
  &\le -\int_0^T\int_\Omega\big(\kappa((w_j)_{k,K})-\kappa(w_{k,K})\big)\na\theta
  \cdot\na(\theta_j-\theta)dxdt \\
  &\phantom{xx}{}
  - \frac{1}{\tau}\int_0^T\int_\Omega((w_j)_K-w_K)\theta(\theta_j-\theta)dxdt \\
  &\phantom{xx}{}
  + \frac{1}{\tau}\int_0^T\int_\Omega((w_j)_K-w_K)\theta_L(\theta_j-\theta)dxdt.
\end{align*}
The regularity $\na\theta\in L^p(0,T;L^p(\Omega))$ for some $p>2$
and the strong convergence
of $\kappa((w_j)_{k,K})\to\kappa(w_{k,K})$ in any $L^r(0,T;L^r(\Omega))$ imply that
$(\kappa((w_j)_{k,K})-\kappa(w_{k,K}))\na\theta\to 0$ strongly 
in $L^2(0,T;L^2(\Omega))$. Hence,
since $\na\theta_j\to\na\theta$ weakly in $L^2(0,T;L^2(\Omega))$, the first
integral on the right-hand side converges to zero. Similarly, in view of 
the $L^\infty$ bounds for $\theta$ and $\theta_L$, the second and third
integrals converge to zero. Since $\kappa((w_j)_{k,K})\ge\kappa_*>0$, 
this shows the claim. 

Hence, we can pass to the limit $j\to\infty$ in the equation
$$
  \int_0^T\langle \pa_t n_j,\phi\rangle dt
  + \int_0^T\int_\Omega\theta_j\na n_j\cdot\na\phi dxdt
  = -\sigma_j\int_0^T\int_\Omega (w_j)_K\na(\theta_j+V_j)\cdot\na\phi dxdt,
$$
where $\langle\cdot,\cdot\rangle$ denotes the dual product on
$H_0^1(\Omega\cup\Gamma_N)$ and $\phi\in L^2(0,T;H_0^1(\Omega\cup\Gamma_N))$, 
to infer that $n$ solves
$$
  \pa_t n - \diver(\theta\na n) 
  = -\sigma\diver(w_K\na(\theta+V))\quad\mbox{in }L^2(0,T;H^{-1}(\Omega\cup\Gamma_N)).
$$
This implies that $n=S(w,\sigma)$. Hence, $S$ is continuous and, by the Aubin lemma,
also compact. 

We prove uniform estimates in $L^\infty(0,T;L^\infty(\Omega))$ 
for all fixed points of $S(\cdot,\sigma)$ which allows to
remove the truncation and which yields uniform estimates in $L^2(0,T;L^2(\Omega))$
needed for the fixed-point theorem.

{\em Step 3: $L^\infty$  bounds for $n$.} 
Let $n$ be a fixed point of $S(\cdot,\sigma)$.
First, observe that the test function $n^-$ in \eqref{3.eq.n} immediately 
implies that $n^-=0$ and $n\ge 0$ in $\Omega$, $t>0$, since $n_K\na n^-=0$
in $\Omega$. To derive an upper bound, we set $u=e^{-\beta t}n$. 
Then $u$ solves the equation
\begin{equation}\label{3.u}
  \pa_t u - \diver(\theta\na u) = \sigma\diver(u_{K_0}\na(\theta+V)) - \beta u,
\end{equation}
since $n_K=\max\{0,\min\{K_0 e^{\beta t},e^{\beta t}u\}\}
=e^{\beta t}\max\{0,\min\{K_0,u\}\} =: e^{\beta t}u_{K_0}$.
Let $L>K$ and define 
$\phi=\kappa(n_{k,K})^{-1} u_{K_0}(u_L-K_0)^+$, where $u_L=\min\{L,u\}$.
This truncation is necessary to obtain $\phi\in L^2(0,T;H_0^1(\Omega\cup\Gamma_N))$.
Furthermore, $\phi(0)=0$ since $u(0)=n_I\le K_0$ in $\Omega$.
We employ the test function $\phi$ in the temperature equation \eqref{3.eq.t}:
\begin{equation}\label{3.phi}
  -\int_\Omega\kappa(n_{k,K})\na\theta\cdot\na\phi dx
  = \frac{1}{\tau}\int_\Omega \frac{n_K}{\kappa(n_{k,K})}(\theta-\theta_L)
  u_{K_0}(u_L-K_0)^+ dx.
\end{equation}
First, we compute the left-hand side:
\begin{align*}
  -\int_\Omega\kappa(n_{k,K})\na\theta\cdot\na\phi dx
  &= -\int_\Omega u_{K_0}\na\theta\cdot\na(u_L-K_0)^+ dx
  - \int_\Omega\na u_{K_0}\cdot\na\theta(u_L-K_0)^+ dx \\
  &\phantom{xx}{}+ \int_\Omega\frac{u_{K_0}}{\kappa(n_{k,K})}(u_L-K_0)^+
  \frac{\pa\kappa}{\pa n}\na n_K\cdot\na\theta dx.
\end{align*}
The second and third integrals vanish since $\na u_{K_0}=0$ and
$\na n_K=0$ on $\{u>K_0\}$. We obtain
$$
  -\int_\Omega\kappa(n_{k,K})\na\theta\cdot\na\phi dx
  = -K_0\int_\Omega \na\theta\cdot\na(u_L-K_0)^+ dx.
$$
Therefore, since $\theta\le M$ and 
$n_K/\kappa(n_{k,K})\le n_{k,K}/\kappa(n_{k,K})\le 1/\kappa_1$ 
(see \eqref{1.kappa}), \eqref{3.phi} becomes
\begin{align}
  -K_0\int_\Omega \na\theta\cdot\na(u_L-K_0)^+ dx
  &= \frac{1}{\tau}\int_\Omega \frac{n_K}{\kappa(n_{k,K})}(\theta-\theta_L)
  u_{K_0}(u_L-K_0)^+ dx \label{3.phi2} \\
  &\le \frac{M}{\tau\kappa_1}\int_\Omega u_{K_0}(u_L-K_0)^+ dx. \nonumber
\end{align}

Next, we use $(u_L-K_0)^+$ as an admissible test function in \eqref{3.u}.
An elementary computation shows that
$$
  F(s) = \int_0^s (\sigma_L-K_0)^+ d\sigma \ge \frac12\big((s_L-K_0)^+\big)^2.
$$
Therefore, since $F(u(0))=F(n_I)=0$,
$$
  \int_0^t\langle \pa_t u,(u_L-K_0)^+\rangle ds
  = \int_\Omega\big(F(u(t))-F(u(0))\big)dx 
  \ge \frac12\int_\Omega\big((u(t)_L-K_0)^+\big)^2 dx,
$$
where $\langle\cdot,\cdot\rangle$ denotes the dual product on
$H_0^1(\Omega\cup\Gamma_N)$. This gives
\begin{align*}
  \frac12\int_\Omega & \big((u(t)_L-K_0)^+\big)^2 dx
  + \int_0^t\int_\Omega \theta|\na(u_L-K_0)^+|^2 dxdt \\
  &\le -\sigma\int_0^t\int_\Omega u_{K_0}\na(\theta+V)\cdot\na(u_L-K_0)^+ dxdt
  - \beta\int_0^t\int_\Omega u(u_L-K_0)^+ dxdt.
\end{align*}
By the Poisson equation \eqref{3.eq.v},
\begin{align*}
  -\int_\Omega u_{K_0}\na V\cdot\na(u_L-K_0)^+ dx
  &= -K_0\int_\Omega \na V\cdot\na(u_L-K_0)^+ dx \\
  &= -\lambda^{-2}K_0\int_\Omega(n_K-C(x))(u_L-K_0)^+ dx
  \le 0,
\end{align*}
since $u>K_0$ is equivalent to $n>K$ and hence, $n_K-C(x)=K-C(x)\ge K_0-C(x)\ge 0$
on $\{u>K_0\}$, using the definition of $K_0$. 
Then, taking into account \eqref{3.phi2}, we find that
\begin{align*}
  \frac12\int_\Omega & \big((u(t)_L-K_0)^+\big)^2 dx
  + m\int_0^t\int_\Omega |\na(u_L-K_0)^+|^2 dxdt \\
  &\le \frac{M}{\tau \kappa_1}\int_0^t\int_\Omega u_{K_0}(u_L-K_0)^+ dx
  - \beta\int_0^t\int_\Omega u(u_L-K_0)^+ dxdt  \\
  &\le \frac{M}{\tau \kappa_1}\int_0^t\int_\Omega u_{K_0}(u_L-K_0)^+ dx
  - \beta\int_0^t\int_\Omega u_{K_0}(u_L-K_0)^+ dxdt  \\
  &= \left(\frac{M}{\tau \kappa_1}-\beta\right)
  \int_0^t\int_\Omega u_{K_0}(u_L-K_0)^+ dxdt
  = 0,
\end{align*}
by the definition of $\beta$. We infer that $(u_L-K_0)^+=0$ for all $L>K_0$.
Letting $L\to\infty$, we obtain $(u-K_0)^+=0$ and thus, $n\le K$
in $\Omega$, $t>0$. As a consequence, $n_K=n$, and any solution to
\eqref{3.eq.n}-\eqref{3.eq.v} solves \eqref{1.eq.n}-\eqref{1.eq.v}.
Furthermore, the $L^\infty$ bounds provide the uniform estimates needed to
apply the Leray-Schauder fixed-point theorem. This proves the existence of solutions
to \eqref{1.eq.n}-\eqref{1.bc}.

{\em Step 4: Positive lower bound for $n$.}
Assume that $\kappa(z)=z$ for all $0\le z\le n_*$. We claim that under
this condition, $n$ possesses a positive lower bound.
In view of the upper bound from Step 3, $(n-k)^-$, where $k=k_0 e^{-\alpha t}$, 
is an admissible test function in \eqref{3.eq.n} yielding
\begin{align}
  \frac12\int_\Omega & (n-k)^-(t)^2 dx
  + m\int_0^t\int_\Omega|\na(n-k)^-|^2 dxdt
  \le -\sigma\int_0^t\int_\Omega n\na\theta\cdot\na(n-k)^- dxdt \label{3.nn} \\
  &{}- \sigma\int_0^t\int_\Omega n\na V\cdot\na(n-k)^- dxdt
  + \alpha\int_0^t\int_\Omega k(n-k)^- dxdt. \nonumber
\end{align}
We write the second integral on the right-hand side as
\begin{align*}
  -\sigma\int_0^t & \int_\Omega (n-k)\na V\cdot\na(n-k)^- dxdt
  -\sigma\int_0^t\int_\Omega k\na V\cdot\na(n-k)^- dxdt \\
  &= -\frac{\sigma}{2}\int_0^t \int_\Omega \na V\cdot\na\big((n-k)^-\big)^2 dxdt
  -\sigma\int_0^t\int_\Omega k\na V\cdot\na(n-k)^- dxdt \\
  &= -\frac{\sigma}{2\lambda^2}\int_0^t\int_\Omega(n-C(x))\big((n-k)^-\big)^2 dxdt
  - \frac{\sigma}{\lambda^2}\int_0^t\int_\Omega (n-C(x))(n-k)^- dxdt \\
  &\le \frac{1}{2\lambda^2}\|C\|_{L^\infty(\Omega)}\int_0^t\int_\Omega 
  \big((n-k)^-\big)^2 dxdt + \frac{1}{\lambda^2}\int_0^t\int_\Omega k[-(n-k)^-] dxdt,
\end{align*}
using the Poisson equation and $n(n-k)^- \le k(n-k)^-$ in $\Omega$. 

In order to estimate the first integral on the right-hand side of \eqref{3.nn}, 
we employ the test function $(n-k)^-$
in \eqref{3.eq.t}. Then, since $\kappa(n)=n$ for $0\le n < k\le k_0\le n_*$,
$$
  \frac{1}{\tau}\int_\Omega n(\theta-\theta_L)(n-k)^- dx
  = -\int_\Omega\kappa(n)\na\theta\cdot\na(n-k)^- dx 
  = -\int_\Omega n\na\theta\cdot\na(n-k)^- dx.
$$
Therefore, \eqref{3.nn} becomes
\begin{align*}
  \frac12\int_\Omega & (n-k)^-(t)^2 dx
  + m\int_0^t\int_\Omega|\na(n-k)^-|^2 dxdt
  \le \frac{\sigma}{\tau}\int_0^t\int_\Omega n(\theta-\theta_L)(n-k)^- dx \\
  &\phantom{xx}{}+ \frac{1}{2\lambda^2}\|C\|_{L^\infty(\Omega)}\int_0^t\int_\Omega 
  \big((n-k)^-\big)^2 dxdt 
  + \left(\frac{1}{\lambda^2}-\alpha\right)\int_0^t\int_\Omega k[-(n-k)^-] dxdt \\
  &\le \frac{1}{2\lambda^2}\|C\|_{L^\infty(\Omega)}\int_0^t\int_\Omega 
  \big((n-k)^-\big)^2 dxdt \\
  &\phantom{xx}{}+ \left(\frac{1}{\tau}\|\theta_L\|_{L^\infty(\Omega)}
  + \frac{1}{\lambda^2}-\alpha\right)\int_0^t\int_\Omega k[-(n-k)^-] dxdt \\
  &= \left(\frac{1}{\tau}\|\theta_L\|_{L^\infty(\Omega)}
  + \frac{1}{\lambda^2}-\alpha\right)\int_0^t\int_\Omega k[-(n-k)^-] dxdt = 0.
\end{align*}
We obtain $(n-k)^-=0$ and hence, $n\ge k$ in $\Omega$, $t>0$.


\section{Proof of Theorem \ref{thm.unique}}\label{sec.unique}

Let $(n_1,\theta_1,V_1)$, $(n_2,\theta_2,V_2)$ be two solutions to
\eqref{1.eq.n}-\eqref{1.eq.v} with the regularity indicated in the theorem.

{\em Step 1: Estimate of $\na(\theta_1-\theta_2)$.}
We employ the test function $\theta_1-\theta_2$ in the difference of the
weak formulations for $\theta_1$, $\theta_2$, respectively:
\begin{align}
  \int_0^t\int_\Omega & \kappa(n_2)|\na(\theta_1-\theta_2)|^2 dxdt
  = -\int_0^t\int_\Omega\big(\kappa(n_1)-\kappa(n_2)\big)
  \na\theta_1\cdot\na(\theta_1-\theta_2)dxdt \label{5.aux1} \\ 
  &\phantom{xx}{}- \frac{1}{\tau}\int_0^t\int_\Omega\big(n_2(\theta_1-\theta_2)
  + (n_1-n_2)(\theta_1-\theta_L)\big)(\theta_1-\theta_2)dxdt. \nonumber
\end{align}
Using the Cauchy-Schwarz, Poincar\'e, and Young inequalities,
the second integral is estimated from above by
\begin{align*}
  c\|n_1 & -n_2\|_{L^2(0,T;L^2(\Omega))} \|\theta_1-\theta_2\|_{L^2(0,T;L^2(\Omega))} \\
  &\le \eps\|\na(\theta_1-\theta_2)\|_{L^2(0,T;L^2(\Omega))}^2
  + c(\eps)\|n_1-n_2\|_{L^2(0,T;L^2(\Omega))}^2.
\end{align*}
where $c(\eps)>0$ depends on $\eps$, the $L^\infty$ bounds for 
$\theta_1$ and $\theta_L$, and the Poincar\'e constant. 
The Lipschitz continuity of $\kappa$ on $[0,K]$ implies that
\begin{align*}
  -\int_0^t\int_\Omega & \big(\kappa(n_1)-\kappa(n_2)\big)
  \na\theta_1\cdot\na(\theta_1-\theta_2)dxdt \\
  &\le c\int_0^t\int_\Omega|n_1-n_2|\,|\na\theta_1|
  \,|\na(\theta_1-\theta_2)|dxdt \\
  &\le c\|n_1-n_2\|_{L^2(0,T;L^2(\Omega))} 
  \|\na(\theta_1-\theta_2)\|_{L^2(0,T;L^2(\Omega))} \\
  &\le \eps\|\na(\theta_1-\theta_2)\|_{L^2(0,T;L^2(\Omega))}^2
  + c(\eps)\|n_1-n_2\|_{L^2(0,T;L^2(\Omega))}^2,
\end{align*}
where $c(\eps)>0$ depends on $\eps$ and the $L^\infty$ norm of $\na\theta_1$.
Since $\kappa(n_2)\ge\kappa_*>0$ for some $\kappa_*>0$, we find from \eqref{5.aux1},
for $\eps\le \kappa_*/4$, that
\begin{equation}\label{5.theta}
  \|\na(\theta_1-\theta_2)\|_{L^2(0,T;L^2(\Omega))}
  \le c(\kappa_*)\|n_1-n_2\|_{L^2(0,T;L^2(\Omega))}.
\end{equation}

{\em Step 2: Estimate of $n_1-n_2$.}
We employ $n_1-n_2$ in the difference of the equations satisfied by
$n_1$ and $n_2$, respectively:
\begin{align}
  \frac12\int_\Omega & (n_1-n_2)(t)^2 dx 
  + \int_0^t\int_\Omega\theta_2|\na(n_1-n_2)|^2 dxdt \label{5.aux4} \\
  &= -\int_0^t\int_\Omega(\theta_1-\theta_2)\na n_1\cdot\na(n_1-n_2)dxdt \nonumber \\
  &\phantom{xx}{}
  + \int_0^t\int_\Omega \big(n_2\na(\theta_1-\theta_2)+(n_1-n_2)\na\theta_1\big)
  \cdot\na(n_1-n_2)dxdt \nonumber \\
  &\phantom{xx}{}
  + \int_0^t\int_\Omega\big(n_2\na(V_1-V_2)+(n_1-n_2)\na V_1\big)\cdot\na(n_1-n_2)dxdt.
  \nonumber 
\end{align}
Applying H\"older's inequality with $p>2$ as in the theorem 
and $1/p+1/q+1/2=1$ to the first integral, we estimate as follows:
\begin{align*}
  -\int_0^t\int_\Omega & (\theta_1-\theta_2)\na n_1\cdot\na(n_1-n_2)dxdt \\
  &\le \|\theta_1-\theta_2\|_{L^2(0,T;L^q(\Omega))}
  \|\na n_1\|_{L^\infty(0,T;L^p(\Omega))}\|\na(n_1-n_2)\|_{L^2(0,T;L^2(\Omega))} \\
  &\le c\|\na(\theta_1-\theta_2)\|_{L^2(0,T;L^2(\Omega))}
  \|\na(n_1-n_2)\|_{L^2(0,T;L^2(\Omega))} \\
  &\le c\|n_1-n_2\|_{L^2(0,T;L^2(\Omega))}\|\na(n_1-n_2)\|_{L^2(0,T;L^2(\Omega))} \\
  &\le \eps\|\na(n_1-n_2)\|_{L^2(0,T;L^2(\Omega))}^2
  + c(\eps)\|n_1-n_2\|_{L^2(0,T;L^2(\Omega))}^2.
\end{align*}
In the second step we have used the Sobolev embedding $H^1(\Omega)\hookrightarrow
L^q(\Omega)$ and the Poincar\'e inquality, 
and the third step follows from \eqref{5.theta}. 

For the second integral in \eqref{5.aux4}, we obtain, using again \eqref{5.theta},
\begin{align*}
  \int_0^t\int_\Omega & \big(n_2\na(\theta_1-\theta_2)+(n_1-n_2)\na\theta_1\big)
  \cdot\na(n_1-n_2)dxdt \\
  &\le \eps\|\na(n_1-n_2)\|_{L^2(0,T;L^2(\Omega))}^2
  + c(\eps)\|\na(\theta_1-\theta_2)\|_{L^2(0,T;L^2(\Omega))}^2 \\
  &\phantom{xx}{}
  + c(\eps)\|n_1-n_2\|_{L^2(0,T;L^2(\Omega))}^2 \\
  &\le \eps\|\na(n_1-n_2)\|_{L^2(0,T;L^2(\Omega))}^2 
  + c(\eps)\|n_1-n_2\|_{L^2(0,T;L^2(\Omega))}^2.
\end{align*}
Finally, for the third integral in \eqref{5.aux4}, we estimate
\begin{align*}
  \int_0^t\int_\Omega & \big(n_2\na(V_1-V_2)+(n_1-n_2)\na V_1\big)\cdot\na(n_1-n_2)dxdt
  \\
  &\le \eps\|\na(n_1-n_2)\|_{L^2(0,T;L^2(\Omega))}^2
  + c(\eps)\|\na(V_1-V_2)\|_{L^2(0,T;L^2(\Omega))}^2 \\
  &\phantom{xx}{}+ \frac12\int_0^t\int_\Omega\na V_1\cdot\na(n_1-n_2)^2 dxdt.
\end{align*}
By the elliptic estimate for the Poisson equation,
\begin{align*}
  \int_0^t\int_\Omega & \big(n_2\na(V_1-V_2)+(n_1-n_2)\na V_1\big)\cdot\na(n_1-n_2)dxdt
  \\
  &\le \eps\|\na(n_1-n_2)\|_{L^2(0,T;L^2(\Omega))}^2
  + c(\eps)\|n_1-n_2\|_{L^2(0,T;L^2(\Omega))}^2 \\
  &\phantom{xx}{}+ \frac{1}{2\lambda^2}\int_0^t\int_\Omega(n_1-C(x))(n_1-n_2)^2 dxdt \\
  &\le \eps\|\na(n_1-n_2)\|_{L^2(0,T;L^2(\Omega))}^2
  + c(\eps)\|n_1-n_2\|_{L^2(0,T;L^2(\Omega))}^2.
\end{align*}
Inserting these estimates in \eqref{5.aux4} and observing that
$\theta_2$ is uniformly bounded from below, i.e.\ $\theta_2\ge m>0$
in $\Omega$, $t>0$, we infer that
\begin{align*}
  \frac12\|(n_1 & -n_2)(t)\|_{L^2(\Omega)}^2
  + m\|\na(n_1-n_2)\|_{L^2(0,T;L^2(\Omega))}^2 \\
  &\le 3\eps\|\na(n_1-n_2)\|_{L^2(0,T;L^2(\Omega))}^2
  + c(\eps)\|n_1-n_2\|_{L^2(0,T;L^2(\Omega))}^2.
\end{align*}
Then, choosing $\eps\le 1/(3m)$, the Gronwall lemma allows us to conclude that
$(n_1-n_2)(t)=0$ in $\Omega$ for $t>0$. This proves the uniqueness of solutions.


\section{Numerical approximation}\label{sec.num}

In this section, we present numerical results for the simplified energy-transport
model with $\kappa(n,\theta)=n\theta$ on the interval $[0,1]$. 
The initial and boundary conditions are
\begin{align*}
  & n_I (x) = C(x)\quad\mbox{for }x\in\Omega, \quad n(0,t) = C(0), 
  \quad n(1,t) = C(1), \\
  & \theta(0,t) = \theta_L(0), \quad \theta(1,t) = \theta_L(1), \quad
  V(0,t) = 0, \quad V(1,t) = U \quad\mbox{for }t>0.
\end{align*}
Equations \eqref{1.eq.n}-\eqref{1.eq.v} are discretized on an equidistant grid
with $N$ grid points $x_i=i\triangle x$, where $\triangle x=1/(N-1)$. The 
time grid points are $t_k=k\triangle t$, where $\triangle t>0$.
We employ central finite differences in space
and the trapezoidal rule in time. Then, with the approximations $n_i^k$,
$\theta_i^k$, and $V_i^k$ of $n(x_i,t_k)$, $\theta(x_i,t_k)$, and $V(x_i,t_k)$,
respectively, the discretized equations become
\begin{align*}
  \frac{1}{\triangle t}(n_i^k-n_i^{k-1})
  &= \frac{1}{(\triangle x)^2}\big((n_{i+1}^{k}\theta_{i+1}^{k} 
  - 2 n_{i}^{k} \theta_{i}^{k} + n_{i-1}^{k} \theta_{i-1}^{k}) \\
  &\phantom{xxxx}{}
  + (n_{i+1}^{k-1} \theta_{i+1}^{k-1} - 2 n_{i}^{k-1} \theta_{i}^{k-1} 
  + n_{i-1}^{k-1} \theta_{i-1}^{k-1})\big) \\
  &\phantom{xx}{}- \frac{1}{2 (\triangle x)^2} \big( 
  (n_{i+1}^k + n_i^k)(V^k_{i+1} - V^k_i) 
  - (n^k_{i} + n^k_{i-1})(V^k_{i} - V^k_{i-1})\big) \\
  &\phantom{xx}{}- \frac{1}{2(\triangle x)^2} \big( 
  (n^{k-1}_{i+1} + n^{k-1}_i)(V^{k-1}_{i+1} - V^{k-1}_i) \\
  &\phantom{xxxx}{}
  - (n^{k-1}_{i} + n^{k-1}_{i-1})(V^{k-1}_{i} - V^{k-1}_{i-1}) \big), \\
  n_i^k - C_i &= \frac{\lambda^2}{(\triangle x)^2}\big(V_{i+1}^k - 2 V_{i}^k 
  + V_{i-1}^k\big), \\
  \frac{n^k_i}{\tau}(\theta^k_i - \theta_{L,i}) 
  &= \frac{\kappa}{2(\triangle x)^2} 
  \big((n_{i+1}^k \theta_{i+1}^k + n_i^k\theta_i^k)(\theta^k_{i+1} - \theta^k_i) 
  - (n^k_{i} \theta^k_{i} + n^k_{i-1} \theta^k_{i-1})(\theta^k_{i} - \theta^k_{i-1})    
  \big).
\end{align*}
Given $(n_i^{k-1},\theta_i^{k-1},V_i^{k-1})$, we find $(n_i^k,\theta_i^k,V_i^k)$
by solving the above nonlinear equations subject to the corresponding (Dirichlet)
boundary conditions using Newton's method. 

We simulate a ballistic diode which is defined by the doping profile
$$
  C(x) = 1 + 0.25\big(\tanh(100x-60) - \tanh(100x-40)\big), \quad x\in[0,1].
$$
The physical parameters are given in Table \ref{table}, and the scaled quantities
are defined by
$$
  \lambda^2 = \frac{\eps_0 \eps_r k_B T_0}{q C_{\rm max} L^2},
  \quad \kappa = \kappa_0 \tau_0 \frac{k_B T_0}{m_n}, \quad t^* =
  \sqrt{\frac{m_n L^2}{k_B T_0}}, \quad \tau = \frac{\tau_0}{t^*}.
$$
For the computations, we choose $N=201$ grid points and the time step size
$\triangle t=1.25\times 10^{-4}$. 

\begin{table}[ht]
\centering
  \begin{tabular}{lll} \hline
  Parameter & Value & Physical meaning  \\ \hline 
  $ k_B$ & $1.3807\times 10^{-23}$\,kg\,m/s$^2$K & Boltzmann constant \\
  $\epsilon_0$ & $8.8542\times 10^{-12}$\,A$^2$\,s$^4$/kg\,m$^3$
  & Vacuum permittivity \\
  $m_0$ & $9.11\times 10^{-31}$\,kg & Electron mass at rest \\
  $q$ & $1.602\times 10^{-19}$\,A\,s& Elementary charge \\
  $C_{\rm max}$ & $10^{24}$\,m$^{-3}$ & Maximum doping concentration\\
  $T_0$ & 300\,K & Device temperature\\
  $L$ & 75\,nm & Device length\\
  $m_n$ & $0.067 \cdot m_0$  & Effective electron mass\\
  $\eps_r$ & $11.7$ & Relative permittivity of GaAs \\
  $\tau_0$ & $0.9\times 10^{-12}$\,s & Momentum relaxation time\\  
  $\lambda^2$ & $3.0\times 10^{-3}$ & Scaled squared Debye length \\
  $\tau$ & 3.126 & Scaled energy relaxation time\\
  $\kappa_0$ & $4.88\times 10^{-2}$ & Heat transfer coefficient\\
\hline
\end{tabular}
\caption{Physical and scaled parameters.}\label{table}
\end{table}

\begin{figure}[ht]
\includegraphics[width=160mm,height=80mm]{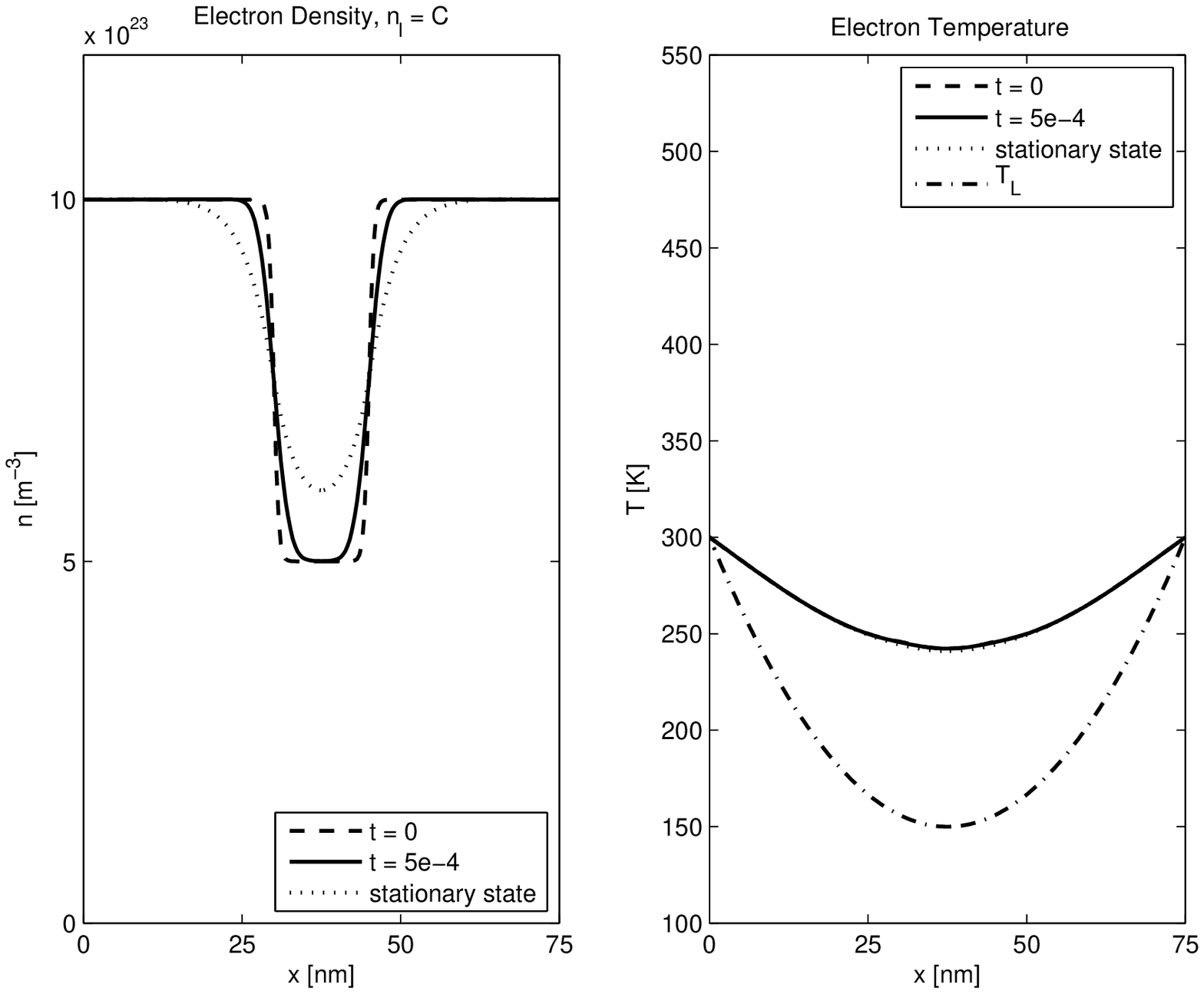}
\includegraphics[width=160mm,height=80mm]{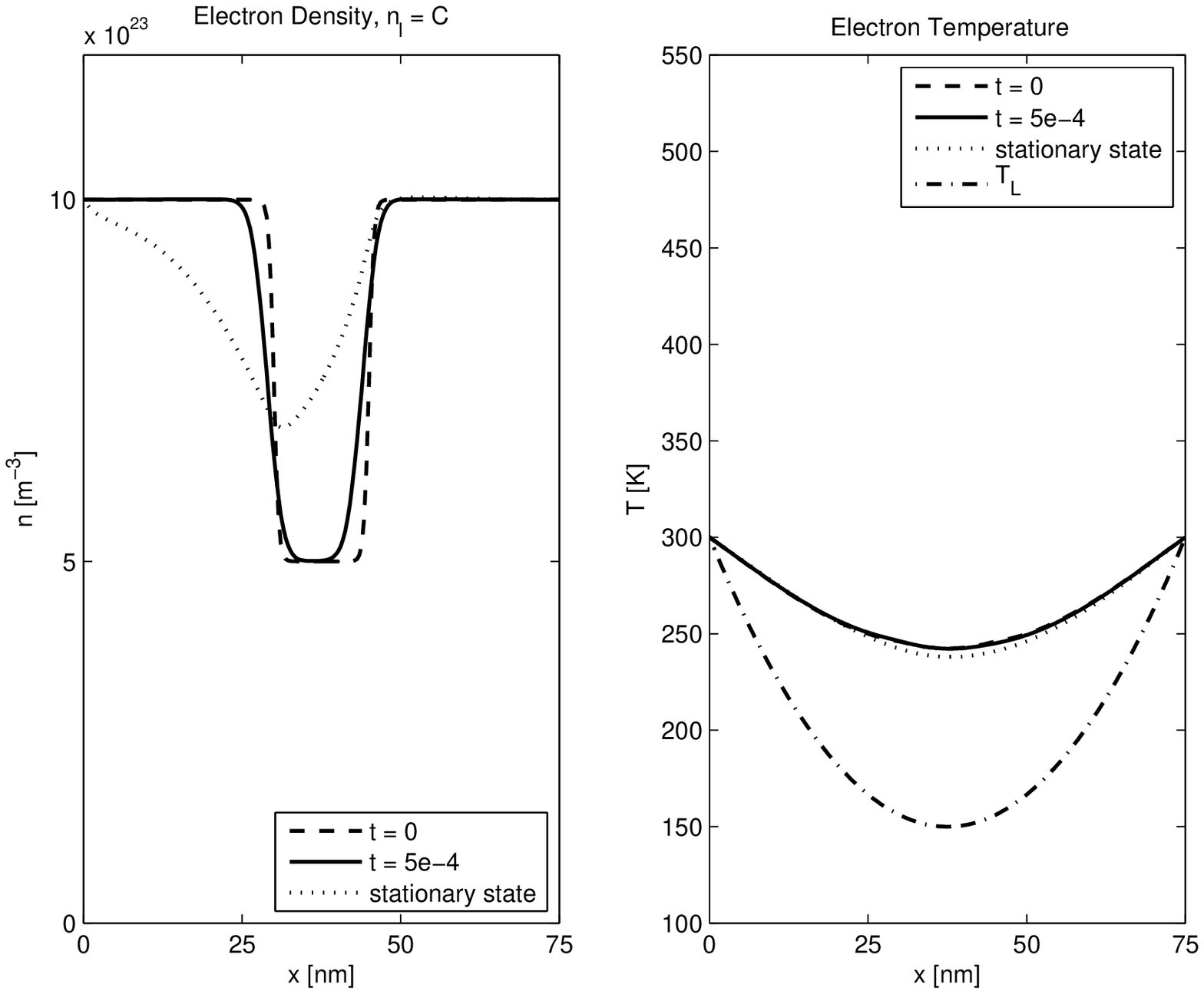}
\caption{Electron density and temperature in the ballistic diode with
cooling lattice temperature at voltages $U=0.2$\,V (top) and 
$U=1.0$\,V (bottom).}\label{fig.c}
\end{figure}

We wish to study the impact of different lattice temperatures. 
First, we choose a lattice temperature which is cooling the interior of the diode, 
i.e.\ $\theta_L(x)=\frac12(x-\frac12)^2+\frac12$.
Figure \ref{fig.c} shows the electron density and electron temperature at various times
for applied voltages $U=0.2$\,V and $U=1.0$\,V, respectively.
In both cases, the electron temperature converges to
its {\em nonhomogeneous} stationary profile as $t\to\infty$. 
Since the profile is convex,
equation \eqref{1.eq.t} implies that the particle temperature is larger than
the lattice temperature. The profile of the electron density follows the
doping profile except for the large applied bias $U=1.0$\,V. In this situation,
the electric force is sufficiently strong to deplete the charge carrier concentration
close to the left boundary point.

\begin{figure}[ht]
\includegraphics[width=160mm,height=80mm]{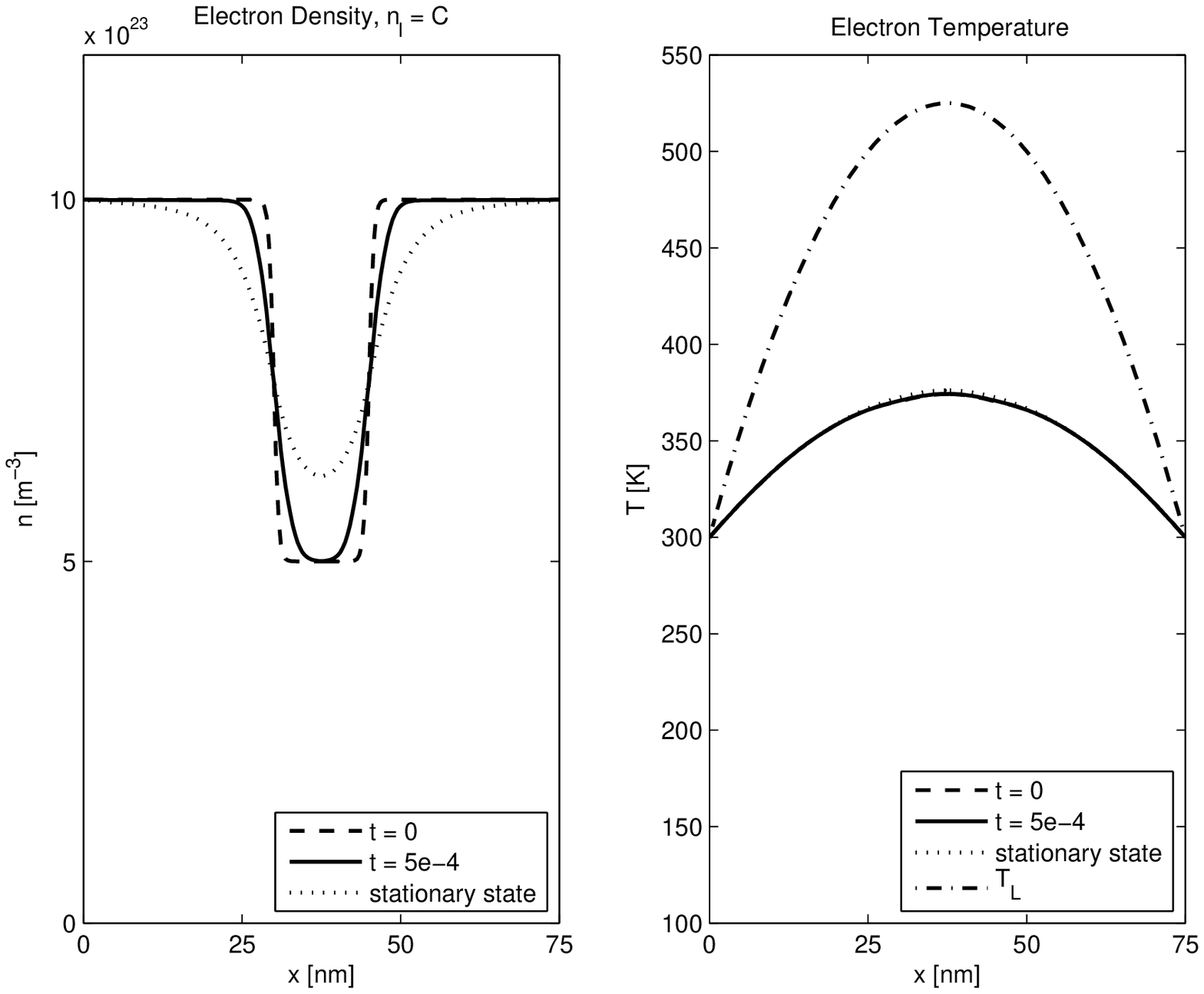}
\includegraphics[width=160mm,height=80mm]{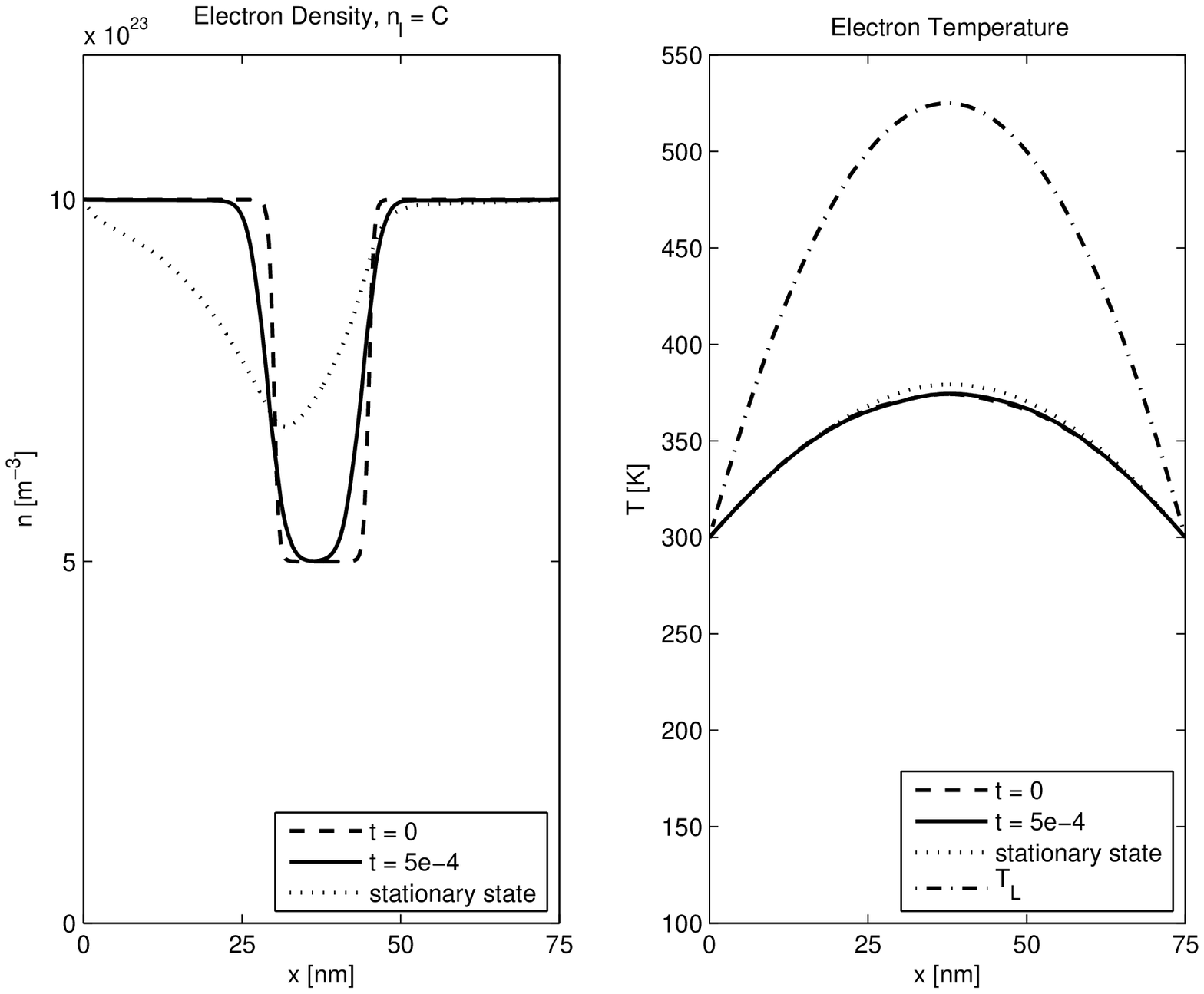}
\caption{Electron density and temperature in the ballistic diode with
heating lattice temperature at voltages $U=0.2$\,V (top) and 
$U=1.0$\,V (bottom).}\label{fig.h}
\end{figure}

Figure \ref{fig.h} illustrates the behavior of the electron density
and electron temperature when the lattice temperature is heating the diode, 
i.e.\ $\theta_L(x)=\frac74-3(x-\frac12)^2$.
Again, the electron temperature converges to a nonhomogeneous steady state, and the 
behavior of the particle density is similar to the case of cooling temperatures.
The current-voltage characteristic is very close to that one with constant
temperature (not presented). It can be seen that only for very large voltages, 
the current density becomes slightly smaller due to an increasing thermal 
energy fraction. This shows that the influence of the temperature equation 
is not very important in a ballistic diode although there are significant 
temperature gradients.

\end{document}